\documentclass{article}
\usepackage{amssymb,verbatim,amsmath,amsthm}

\newtheorem{thm}{Theorem}[section]
\newtheorem{prop}[thm]{Proposition}
\newtheorem{conj}[thm]{Conjecture}
\newtheorem{lem}[thm]{Lemma}

\theoremstyle{definition}
\newtheorem{defn}[thm]{Definition}

\theoremstyle{remark}

\title{Extremal metrics and K-stability}

\author{G\'abor Sz\'ekelyhidi}
\date{}

\begin{document}

\maketitle
\begin{abstract}

We propose an algebraic geometric stability criterion for
a polarised variety to admit an extremal K\"ahler metric. This generalises
conjectures by Yau, Tian and Donaldson which relate to the case of
K\"ahler-Einstein and constant scalar curvature metrics. 
We give a result in geometric invariant theory
that motivates this conjecture, and an
example computation that supports it.

\end{abstract}

\section{Introduction} \label{sec:intr}

Much progress has been made recently in understanding the relation
between stability, and existence of special metrics on smooth algebraic
varieties. Such a relation was originally conjectured by Yau for
K\"ahler-Einstein metrics. 
One of the problems is to find the appropriate notion of
stability.  An example is K-stability introduced by Tian~\cite{Tian97}
(see also Donaldson~\cite{Don02}) which is conjectured to be equivalent
to the existence of a K\"ahler-Einstein metric, or more generally, a
K\"ahler metric of constant scalar curvature. In this paper, we
define a modification of K-stability and formulate a conjecture
relating it to the existence of extremal metrics in the sense of
Calabi~\cite{Cal82}. These metrics are defined to be critical points of
the $L^2$ norm of the scalar curvature defined on K\"ahler metrics in
a fixed cohomology class. See also Mabuchi~\cite{Mab04_1} for the relation
between stability and extremal metrics.

The definition of K-stability in \cite{Don02} involves considering
test-configurations of the
variety which are degenerations
into possibly singular and non-reduced schemes. For every
test-configuration there is an
induced $\mathbf{C}^*$-action on the central fibre, and the stability
condition is that the corresponding generalised Futaki invariant is
non-negative (using the convention in Ross-Thomas~\cite{RT04} for the
Futaki invariant). 

Our modification essentially consists of working ``orthogonal'' to a
maximal torus of automorphisms of the variety. This means that we only
consider test-con\-fi\-gu\-ra\-tions which commute with this torus and the
$\mathbf{C}^*$-action induced on the central fibre is orthogonal to the
torus with respect to an inner product to be defined. Since this
orthogonality condition is not very natural, it is more convenient to
instead modify the Futaki invariant as follows. The maximal torus
contains a distinguished $\mathbf{C}^*$-action $\chi$, generated by the
extremal vector field (see Futaki and Mabuchi~\cite{FM95}). This induces
an action on the central fibre which we also denote by $\chi$. 
Denote by $\alpha$ the $\mathbf{C}^*$-action on the central fibre
induced by the test-configuration and write $F(\alpha)$ for the
Futaki invariant. The modified Futaki invariant is then
\begin{equation}
  \label{eq:modfutaki}
  F_\chi(\alpha) = F(\alpha) - \langle\alpha,\chi\rangle.
\end{equation}

\noindent The variety is K-stable relative to the maximal torus if
$F_\chi(\alpha)$ is non-negative for all test-configurations commuting
with the torus, and zero only for test-configurations arising from a
$\mathbf{C}^*$-action on the variety. 

\begin{conj}\label{mainconj} A polarised variety admits an extremal
  metric in the class of the polarisation, if
and only if it is K-stable relative to a maximal torus.
\end{conj}

In the next section we will give the definitions of the Futaki invariant
and the inner product. In section~\ref{sec:moment} we motivate
Conjecture \ref{mainconj}
by a result in the finite dimensional framework of moment
maps and stability. Finally in section~\ref{sec:ex} we give an example
to support Conjecture \ref{mainconj}.

\subsection*{Acknowledgements} I would like to thank my PhD supervisor
Simon Donaldson for introducing me to this problem, and 
Richard Thomas for helpful conversations and for his many
suggestions improving earlier versions of this paper. I would also like to
thank the referee for his comments, and EPSRC for financial support.

\section{Basic definitions} \label{sec:defn}

We first recall the definition of the generalized Futaki invariant from
Donaldson~\cite{Don02}.
Let $V$ be a polarised scheme of dimension $n$, with a very ample line
bundle $\mathcal{L}$.  Let
$\alpha$ be a $\mathbf{C}^*$-action on $V$ with a lifting to
$\mathcal{L}$.  This induces
a $\mathbf{C}^*$-action on the vector space of sections
$H^0(V,\mathcal{L}^k)$ for
all integers $k\geqslant 1$. Let $d_k$ be the dimension of
$H^0(V,\mathcal{L}^k)$,
and denote the infinitesimal generator of the
action by $A_k$. Denote by $w_k(\alpha)$ the weight of the action on the top
exterior power of $H^0(V,\mathcal{L}^k)$. This is the same as the trace
$\mathrm{Tr}(A_k)$. Then $d_k$ and $w_k(\alpha)$ are polynomials
in $k$ of degree $n$ and $n+1$ respectively for $k$ sufficiently large,
so we can write  
\begin{eqnarray*}
d_k&=&c_0k^n + c_1k^{n-1} + O(k^{n-2})\\
w_k(\alpha)=\mathrm{Tr}(A_k)&=& a_0k^{n+1} + a_1k^n + O(k^{n-1}).
\end{eqnarray*}

The \emph{Futaki invariant} is defined to be
$F(\alpha)=a_0c_1-a_1c_0$. 
The choice of lifting of $\alpha$ to the line
bundle is not unique, however $A_k$ is defined up to addition of a
scalar matrix. In fact if we embed $V$ into $\mathbf{P}^{d_1-1}$ using
$\mathcal{L}$, then lifting $\alpha$ is equivalent to
giving a $\mathbf{C}^*$-action on $\mathbf{C}^{d_1}$ which induces
$\alpha$ on $V$ in
$\mathbf{P}^{d_1-1}$. Since
the embedding by sections of a line bundle is not
contained in any hyperplane, two such $\mathbf{C}^*$-actions differ by
an action that acts trivially on $\mathbf{P}^{d_1-1}$ ie. one with a
constant weight, say $\lambda$. We obtain that for
another lifting, the sequence of matrices $A^\prime_k$ are related to
the $A_k$ by 
\[ A^\prime_k = A_k + k\lambda I,\]
where $I$ is the identity matrix. A simple computation shows
that $F(\alpha)$ is independent of the lifting of $\alpha$ to
$\mathcal{L}$. 

We are now going to define an inner product  on 
$\mathbf{C}^*$-actions. Since $\mathbf{C}^*$-actions do not naturally form a
vector space, this is not really an inner product, but for
smooth varieties we will see that it is the
restriction of an inner product on a space of holomorphic vector fields 
to the set of $\mathbf{C}^*$-actions.

Let $\alpha$ and $\beta$ be two
$\mathbf{C}^*$-actions on $V$ 
with liftings to $\mathcal{L}$. If we denote the infinitesimal generators of
the actions on $H^0(V,\mathcal{L}^k)$ by $A_k, B_k$, then
$\mathrm{Tr}(A_kB_k)$ is a polynomial of degree $n+2$ in $k$. We define the
inner product $\langle\alpha,\beta\rangle$ to be the leading coefficient in
\begin{eqnarray*}
  &&\mathrm{Tr}\left[
  \left(A_k-\frac{\mathrm{Tr}(A_k)}{d_k}I\right)\left(B_k-
  \frac{\mathrm{Tr}(B_k)}{d_k}I\right)
\right] = \\ &&\qquad = \mathrm{Tr}(A_kB_k) -
\frac{w_k(\alpha)w_k(\beta)}{d_k} =  
\langle\alpha,\beta\rangle k^{n+2} + O(k^{n+1})\qquad \mbox{for }k\gg1.
\end{eqnarray*}
Again, this does not depend on the particular liftings of $\alpha$ and
$\beta$ to the line bundle since we are normalizing each $A_k$ and $B_k$
to have trace zero.

Before we proceed, it is worth looking at the case when the variety is
smooth. In this case we can consider the algebra of holomorphic vector
fields on $V$ which lift to $\mathcal{L}$. This is the Lie algebra
of a group of holomorphic automorphisms of $V$. Inside this group, let
$G$ be the complexification of a maximal compact subgroup $K$. Let
$\mathfrak{g},\mathfrak{k}$ be the Lie algebras of $G,K$. Denoting by
$\mathfrak{k}_\mathbf{Q}$ the elements in $\mathfrak{k}$
which generate circle subgroups, our inner product on
$\mathbf{C}^*$-actions gives an inner
product on $\mathfrak{k}_\mathbf{Q}$. Since this is a dense subalgebra of
$\mathfrak{k}$, the inner product extends to 
$\mathfrak{k}$ by continuity. We further extend
this inner product to $\mathfrak{g}$ by complexification
and compute it differential geometrically. This is analogous
to the computation in Donaldson~\cite{Don02} (section 2.2) showing the
relation between the Futaki invariant as defined here, and the original
differential-geometric definition of Futaki~\cite{Fut83}. Note that
$\mathfrak{g}$ is a space of holomorphic vector fields on $V$ which lift
to $\mathcal{L}$. Let $v,w$ be
two holomorphic vector fields on $V$, with liftings $\hat{v}, \hat{w}$
to $\mathcal{L}$. Let $\omega$ be a K\"ahler metric on $V$ in the class
$2\pi c_1(\mathcal{L})$, induced by a choice of Hermitian metric on
$\mathcal{L}$.
We can then write  
\[ \hat{v} = \overline{v} + if\underline{t}, \qquad
\hat{w}=\overline{w}+ig\underline{t},\]

\noindent where $\overline{v}$ (respectively $\overline{w}$) is the horizontal
lift of $v$ (respectively $w$), $\underline{t}$ is the canonical vector field on the
total space of $\mathcal{L}$, defined by the action of scalar
multiplication, and $f,g$ are smooth functions on $V$. As
in \cite{Don02} we have that 
\[ \overline{\partial}f = -(i_v(\omega))^{0,1}, \qquad
\overline{\partial}g = -(i_w(\omega))^{0,1}, \]

\noindent so in particular $f$ and $g$ are defined up to an additive
constant, and we can normalize them to have zero integral over $V$. We
would like to show that
\[ \langle v,w\rangle = \int_V fg\frac{\omega^n}{n!}, \]

\noindent where we have assumed that $f,g$ have zero integral over $V$.
Making use of the identity $\langle v,w\rangle = \frac{1}{2}(\langle
v+w,v+w\rangle - \langle v,v\rangle - \langle w,w\rangle)$ it is enough
to show this when $v=w$. Furthermore, we can assume that $v$ generates a
circle action since $\mathfrak{k}_\mathbf{Q}$ is dense in
$\mathfrak{k}$.  

We can find the leading coefficients of $d_k,\mathrm{Tr}(A_k),
\mathrm{Tr}(A_kA_k)$ for this circle action using the
equivariant Riemann-Roch formula, in the 
same way as was done in~\cite{Don02}. We find that these leading
coefficients are given by
\[ \int_V\frac{\omega^n}{n!},\quad \int_V f\frac{\omega^n}{n!},\quad
\int_V f^2\frac{\omega^n}{n!},\]

\noindent respectively. If we normalize $f$ to have zero integral over
$V$, then we obtain the formula for the inner product that we were
after. This inner product on holomorphic vector fields has also been
defined by Futaki and Mabuchi in~\cite{FM95},
where it is shown that it
only depends on the K\"ahler class, not the specific representative
chosen. This can also been seen from the fact that it can be defined
algebro-geometrically, just like the Futaki invariant.

We next recall the notion of a test-configuration from~\cite{Don02} and
introduce the modification that we need.

\begin{defn} 
  A \emph{test-configuration for $(V,L)$ of exponent $r$}
  consists of a $\mathbf{C}^*$-equivariant flat family of schemes
  $\pi:\mathcal{V}\to\mathbf{C}$ (where $\mathbf{C}^*$ acts on
  $\mathbf{C}$ by multiplication) and a $\mathbf{C}^*$-equivariant ample
  line bundle $\mathcal{L}$ over $\mathcal{V}$.  We require that the
  fibres $(\mathcal{V}_t,\mathcal{L}|_{\mathcal{V}_t})$ are isomorphic
  to $(V,L^r)$ for $t\not=0$, where $\mathcal{V}_t=\pi^{-1}(t)$. The
  test-configuration is called a \emph{product configuration} if
  $\mathcal{V} = V\times\mathbf{C}$. 

  We say that the test-configuration is \emph{compatible with a torus $T$
  of automorphisms of $(V,L)$}, if there is a torus action on
  $(\mathcal{V},\mathcal{L})$ which preserves the fibres of
  $\pi:\mathcal{V}\to\mathbf{C}$, commutes with the
  $\mathbf{C}^*$-action, and restricts to $T$ on
  $(\mathcal{V}_t,\mathcal{L}|_{\mathcal{V}_t})$ for $t\not=0$. 
\end{defn}
 
With these preliminaries we can state the main definition.

\begin{defn}\label{def:kstable}
  A polarised variety $(V,L)$ is \emph{K-stable relative to a
  maximal torus of automorphisms} if 
  $F_{\tilde{\chi}}(\tilde{\alpha})\geq 0$ for all
test-configurations compatible with the torus, and equality holds only if
the test-configuration is a product configuration. Here we denote by
$\tilde{\alpha}$
and $\tilde{\chi}$ the $\mathbf{C}^*$-actions induced on the central fibre of
the test-configuration ($\tilde{\chi}$ being induced by the extremal
$\mathbf{C}^*$-action $\chi$ in the chosen maximal torus) and
$F_{\tilde{\chi}}(\tilde{\alpha})$ is defined as in (\ref{eq:modfutaki}).
\end{defn}

\section{Moment map and stability} \label{sec:moment}

The aim of this section is to describe a result in the finite
dimensional picture of moment maps and stability, which motivates
Conjecture \ref{mainconj}. First we introduce the necessary notation. Let
$X$ be a finite dimensional K\"ahler variety. When formally applying the results of this
section to the problem of extremal metrics, $X$ will be an infinite
dimensional space of complex 
structures on a variety $V$. The details will be discussed at the end of
this section. Denote the K\"ahler form by $\omega$ and let
$\mathcal{L}$ be a line bundle over $X$ with first Chern class
represented by $\omega$. Suppose a compact connected group $K$ acts on $X$ by
holomorphic transformations, preserving $\omega$, and there is a moment
map for the action
\[ \mu : X \rightarrow \mathfrak{k}^*, \]

\noindent where $\mathfrak{k}^*$ is the dual of the Lie algebra of $K$. 
This allows us to define an action of $\mathfrak{k}$ on sections of
$\mathcal{L}$ as follows. Choose a Hermitian metric on $\mathcal{L}$
such that the corresponding unitary connection has curvature
form given by $-2\pi i\omega$. If $\xi\in\mathfrak{k}$ induces a vector
field $v$ on $X$, and
$f:X\to\mathbf{R}$ is the corresponding Hamiltonian function given by
the composition 
\[ X\xrightarrow{\mu}\mathfrak{k}^*\xrightarrow{\xi}\mathbf{R},\]
then $\xi$ acts on $\mathcal{L}$ via \( \bar{v} + 2\pi f\underline{t}, \)
where $\bar{v}$ is the horizontal lift of $v$, and $\underline{t}$ is the
vertical vector field generating the 
$U(1)$-action on the fibres
(see Donaldson and Kronheimer~\cite{DK90}, section 6).

Suppose that there is a complexification $G$ of the group $K$, with Lie
algebra $\mathfrak{g}=\mathfrak{k}\oplus i\mathfrak{k}$. The action of
$\mathfrak{k}$ on $X$ and $\mathcal{L}$ extends to actions of
$\mathfrak{g}$ by complexification. We will assume that
this infinitesimal action gives rise to an action of $G$ on the pair
$(X,\mathcal{L})$. This is the situation studied in geometric invariant
theory, and the main definition we need is the following.

\begin{defn} 
  A point $x\in X$ is \emph{stable} for the
  action of $G$ on $(X,\mathcal{L})$, if for a choice of lifting
  $\tilde{x}\in\mathcal{L}$ of $x$, the set $G\tilde{x}$ is closed in
  $\mathcal{L}$  
\end{defn}

This notion is what some authors call polystability, and usually
stability requires in addition that the point in question has a
zero-dimensional stabilizer in $G$. The relation between the moment map
and stability is given by the following well-known result (see for
example Mumford-Fogarty-Kirwan~\cite{MFK94}).

\begin{prop}\label{prop:kempf} 
  A point $x\in X$ is stable for the action of $G$ on
  $(X,\mathcal{L})$ if and only if there is an element $g\in G$ such
  that $\mu(g\cdot x)=0$.
\end{prop}
We would like to extend this characterisation of the $G$-orbits of
zeros of the moment map to $G$-orbits of critical points of the norm
squared of the moment map. More
precisely, suppose there is a non-degenerate inner product
$\langle\cdot,\cdot\rangle$ on $\mathfrak{k}$, invariant under the
adjoint action. We also assume that on the Lie algebra $\mathfrak{t}$ of
a maximal torus in $K$, the inner product is rational when restricted to
the kernel of the exponential map. 
Since any two maximal tori are conjugate, it follows that the
inner product is rational in this sense on the Lie algebra of any other
torus in $K$ as well. Using the inner product, identify $\mathfrak{k}$ with
$\mathfrak{k}^*$ and from now on consider the moment map as a map 
from $X$ to $\mathfrak{k}$. The norm squared of the moment map then
defines a function $f:X\rightarrow \mathbf{R}$,
\[ f(x) = \Vert\mu(x)\Vert^2.\]
We will show that the $G$-orbits of critical points of $f$ are
characterized by stability with respect to the action of a certain subgroup of
$G$. 

First of all, by differentiating $f$, we find that $x$ is a
critical point if and only if
the vector field induced by $\mu(x)$ vanishes at $x$. In particular the
minima of the functional are given by $x$ with $\mu(x)=0$ and the other
critical points have nontrivial isotropy groups containing the group
generated by $\mu(x)$. This is a circle subgroup by the following lemma. 

\begin{lem}\label{lem:rat} If $x\in X$ is a non-minimal critical point
  of $f$, then
  $\mu(x)$ generates a circle subgroup of $K$.
\end{lem}
\begin{proof}
  Let $\beta=\mu(x)$ and denote
  by $T$ the closure of the subgroup of $K$ generated by $\beta$. This
  is a compact connected Abelian Lie group, hence it is a torus. 
  Letting $\mathfrak{t}$ be the Lie algebra of $T$, the moment map
  $\mu_T$ for the action of $T$ on $X$ is given by the composition of
  $\mu$ with the orthogonal projection from $\mathfrak{k}$ to
  $\mathfrak{t}$. Since by definition, $\beta\in\mathfrak{t}$, we have
  that $\mu(x)=\mu_T(x)$. Let $v_1,\ldots,v_k$ be an integral basis for
  the kernel of the exponential map from $\mathfrak{t}$ to $T$. Because
  of the rationality assumption on the inner product, what we need
  to show is that $\langle\mu_T(x),v_i\rangle$ is rational for all
  $i$. Since $f_i=\langle\mu_T,v_i\rangle$ is the Hamiltonian function
  for the vector field induced by $v_i$, we know that $v_i$ acts on the
  fibre $\mathcal{L}_x$ via $2\pi f_i(x)\underline{t}$. Since $\exp(v_i)=1$,
  we find that $f_i(x)$ must be an integer. 
\end{proof}

Now we define the subgroups of $G$ which will feature in the
stability condition. For a torus $T$ in $G$ with Lie algebra
$\mathfrak{t}$,  define two subalgebras of
$\mathfrak{g}$:
\begin{eqnarray*}
  \mathfrak{g}_T &:=& \{\alpha\in\mathfrak{g}\,\vert\,
  [\alpha,\beta]=0\quad\text{for all }\beta\in\mathfrak{t}\} \\
  \mathfrak{g}_{T^\perp} &:=& \{\alpha\in\mathfrak{g}_T\,\vert\,
  \langle\alpha,\beta\rangle = 0\quad\text{for all }
  \beta\in\mathfrak{t}\} \subset\mathfrak{g}_T.
\end{eqnarray*}

Denote the corresponding connected subgroups by $G_T$ and $G_{T^\perp}$.
Then $G_T$ is the identity component of the centraliser of $T$ and $G_{T^\perp}$ is a subgroup isomorphic to the quotient
of $G_T$ by $T$. 
Working on the level of the compact subgroup $K$, if
$\mathfrak{t}\subset\mathfrak{k}$, then the same formulae
define Lie algebras $\mathfrak{k}_T, \mathfrak{k}_{T^\perp}$ and
subgroups $K_T, K_{T^\perp}$ of $K$, such that
\begin{eqnarray*}
\mathfrak{k}_T = \mathfrak{k}\cap\mathfrak{g}_T,&\quad
\mathfrak{k}_{T^\perp} =
\mathfrak{k}\cap\mathfrak{g}_{T^\perp}\\
K_T = K\cap G_T,&\quad K_{T^\perp} = K\cap G_{T^\perp}.
\end{eqnarray*}

We can now write down the stability condition that we need.

\begin{defn} 
  Let $T$ be a torus in $G$ fixing $x$. We say that $x$ is \emph{stable
  relative to $T$}, if it is stable for the action of $G_{T^\perp}$ on
  $(X,\mathcal{L})$. 
\end{defn}

The main result of this section is the following. 

\begin{thm}\label{thm:stab}
  A point $x$ in $X$ is in the $G$-orbit of a non-minimal critical point of
  $f$, if and only if it is stable relative to a maximal torus
  which fixes it. 
\end{thm}

Before giving the proof, consider the effect of varying the maximal compact
subgroup of $G$. If we replace $K$ by a conjugate
$gKg^{-1}$ for some $g\in G$ and we replace $\omega$ by $(g^{-1})^*\omega$,
then we obtain a new compact group acting by sympectomorphisms. The associated
moment map $\mu_g$ is related to $\mu$ by
\begin{equation} \label{eq:moment}
\mu_g(gx) = \mathrm{ad}_g\mu(x) \in \mathrm{ad}_g\mathfrak{k},
\end{equation}

\noindent where we identify the Lie algebra of $gKg^{-1}$ with
$\mathrm{ad}_g\mathfrak{k}\subset\mathfrak{g}$. Using the inner
product on $\mathrm{ad}_g\mathfrak{k}$ induced from the bilinear form on
$\mathfrak{g}$, define the function $f_g(x) = \Vert\mu_g\Vert^2$. This
satisfies $f_g(gx)=f(x)$ by (\ref{eq:moment}) and the
$\mathrm{ad}$-invariance of the bilinear form, so in particular the
critical points of $f_g$
are obtained by applying $g$ to the critical points of $f$.

\begin{proof}[of Theorem~\ref{thm:stab}]

  Suppose first that $x$ is in the $G$-orbit of a critical point of $f$.
  By replacing $K$ with a conjugate if necessary, we can assume that $x$
  itself is a critical point, so $\mu(x)$ fixes $x$. By
  Lemma~\ref{lem:rat} we obtain a circle action fixing $x$, generated by
  $\beta=\mu(x)$. Choose a maximal torus $T$ fixing $x$, containing this
  circle. Since the moment map $\mu_{T^\perp}$ for the action of
  $K_{T^\perp}$ on $X$ is the composition of $\mu$ with the
  orthogonal projection from $\mathfrak{k}$ to
  $\mathfrak{k}_{T^\perp}$, we have that $\mu_{T^\perp}(x)=0$.
  By Proposition~\ref{prop:kempf} this implies that $x$ is stable for
  the action of $G_{T^\perp}$. 

  Conversely, suppose $x$ is stable for the action of $G_{T^\perp}$ for
  a maximal torus $T$ which fixes $x$.  Choose a maximal compact
  subgroup $K$ of $G$ containing $T$. Then $K_{T^\perp}$ is a maximal
  compact subgroup of $G_{T^\perp}$ and using the assumption on $x$,
  Proposition~\ref{prop:kempf} implies that $y=gx$ is in the kernel of
  the corresponding moment map $\mu_T$, for some $g\in G_{T^\perp}$.
  Then, for the moment map corresponding to $K$, $\mu(y)$ is contained
  in $\mathfrak{t}$ therefore fixes $y$. This means that $y$ is a
  critical point of $f$.  
\end{proof}

We will now explain how the formula \ref{eq:modfutaki} arises. Since
$G_x$ fixes $x$, the action on the fibre defines a map
$G_x\to\mathbf{C}^*$. The derivative at the identity gives a linear map
$\mathfrak{g}_x\to \mathbf{C}$ which we denote by $-F_x$ in order to
match with the sign of the Futaki invariant.  We 
say that $-F_x(\alpha)$ is the weight of the action of $\alpha$ on
$\mathcal{L}_x$. According to the Hilbert-Mumford numerical criterion
for stability (see~\cite{MFK94}), we have the following necessary and
sufficient condition for a point $x$ to be stable: for all one-parameter
subgroups $t\mapsto\exp(t\alpha)$ in $G_{T^\perp}$, the weight on
the central fibre $\mathcal{L}_{x_0}$ is negative, or equal to zero if
$\exp(t\alpha)$ fixes $x$. Here $x_0$ is defined to be $\lim_{t\to
0}\exp(t\alpha)x$. In other words, the condition is that 
\[F_{x_0}(\alpha) \geqslant 0 ,\]

\noindent with equality if and only if $x$ is fixed by the one-parameter
subgroup. 

It is in\-con\-venient to restrict atten\-tion to one-parameter subgroups in
$G_{T^\perp}$ because the orthogonality condition is not a natural
one for test-configurations. We would therefore like to be able to
consider one-parameter subgroups in $G_T$ and adapt the numerical
criterion. For a one-parameter subgroup in $G_T$ generated by
$\alpha\in\mathfrak{g}_T$ we consider the one-parameter subgroup
in $G_{T^\perp}$ generated by the orthogonal projection of $\alpha$
onto $\mathfrak{g}_{T^\perp}$, which we denote by
$\overline{\alpha}$. We have 
\[ \overline{\alpha} = \alpha - \sum_{i=1}^k\frac{\langle\alpha,\beta_i\rangle}{
\langle\beta_i,\beta_i\rangle}\beta_i,\]
where $\beta_1,\ldots,\beta_k$ is an orthonormal basis for
$\mathfrak{t}$. 
Since $[\alpha,\mathfrak{t}]=0$ and $x$ is fixed by $T$, the central fibre
for the two one-parameter groups generated by $\alpha$ and
$\overline{\alpha}$ is the same, the only difference is the weight of
the action on this fibre. Since $F_{x_0}$ is linear, we obtain
\[ F_{x_0}(\overline{\alpha}) = F_{x_0}(\alpha) -
\sum_{i=1}^k\frac{\langle\alpha, 
\beta_i\rangle}{\langle\beta_i,\beta_i\rangle}F_{x_0}(\beta_i). \]

The extremal vector field $\chi$ is defined to be the
element in $\mathfrak{t}$ dual to the functional $F_x$, restricted to
$\mathfrak{t}$, under the inner product. In other words,
$F_x(\alpha)=\langle\alpha,\chi\rangle$ for all
$\alpha\in\mathfrak{t}$. If we now choose the orthonormal basis
$\beta_i$ such that $\beta_1=\chi/\Vert\chi\Vert$, then the previous
formula reduces 
to 
\[ F_{x_0}(\overline{\alpha}) =
F_{x_0}(\alpha)-\langle\alpha,\chi\rangle. \]

If we define this expression to be $F_{x_0,\chi}(\alpha)$, then the
stability condition is 
equivalent to $F_{x_0,\chi}(\alpha)\geqslant 0$ 
for all one-parameter subgroups generated by
$\alpha\in\mathfrak{g}_{T}$ with equality only if the one-parameter
subgroup fixes $x$. We therefore obtain the following

\begin{thm}
  A point $x\in X$ is in the $G$-orbit of a critical point of $f$, if
  and only if for each one-parameter subgroup of $G$ generated by an
  element $\alpha\in\mathfrak{g}_T$ we have
  \[ F_{x_0,\chi}(\alpha)\geqslant 0, \]
  with equality only if $\alpha$ fixes $x$. Here $T$ is a maximal torus
  fixing $x$ and $\chi$ is the corresponding extremal vector field. 
\end{thm}

To conclude this section we explain why this result motivates
Conjecture~\ref{mainconj}. The main idea is the infinite dimensional picture
described in Donaldson~\cite{Don97}, in which the scalar curvature
arises as a moment map. We start with a
symplectic manifold $M$ with symplectic form $\omega$ and assume for
simplicity that $H^1(M)=0$. The space $X$ is the space of
integrable complex structures on $M$ compatible with $\omega$. Then $X$
is endowed with a natural K\"ahler
metric. Together with $\omega$, the points of $X$ define metrics on $M$,
so $X$ can also be thought of as a space of K\"ahler metrics on $M$. 
The group $K$ is the identity
component of the group of symplectomorphisms of $M$. This acts on $X$,
preserving the symplectic form. The Lie algebra
$\mathfrak{k}$ of
$K$ can be identified with $C^\infty_0(M,\mathbf{R})$, the space of
smooth real valued functions on $M$ with zero integral, using the
Hamiltonian construction (we use the
condition $H^1(M)=0$ here). The dual $\mathfrak{k}^*$ can also be
identified with $C^\infty_0(M,\mathbf{R})$ using the $L^2$ pairing. Then
a moment map for the action of $K$ of $X$ is given by mapping a point in
$X$ to the scalar curvature function of the corresponding metric on $M$. 

The complexification $\mathfrak{g}$ of $\mathfrak{k}$ is
$C^\infty_0(M,\mathbf{C})$ with the $L^2$ product. The corresponding
group $G$ does not exist, but we can consider a foliation of $X$
generated by the action of $\mathfrak{g}$ on $X$, whose leaves would be
the orbits of $G$. As explained in \cite{Don97}, these leaves correspond
to metrics on $M$ in a fixed K\"ahler class. The problem of finding
critical points of the norm squared of the moment map in a ``$G$-orbit''
is therefore the problem of finding extremal metrics in a K\"ahler
class. By the finite dimensional result in this section we expect that
an analogous stability condition will characterize the K\"ahler classes
which contain an extremal metric. 

Fixing an element $J\in X$, elements of the Lie algebra $\mathfrak{g}$
give rise to vector fields on 
$(M,\omega,J)$. If we identify $\mathfrak{g}$ with $C^\infty_0(M,\mathbf{C})$, then
given an element $f+ig\in\mathfrak{g}$ with both $f,g$ real valued, the
corresponding vector field is $V_f + JV_g$. Here $V_f,V_g$ are the
Hamiltonian vector fields corresponding to $f,g$. The $L^2$ inner
product on $\mathfrak{g}$ is therefore an inner product on a space of
vector fields (not necessarily holomorphic) on $M$. What we did in
section~\ref{sec:defn} was to shift attention to the central
fibre of a test-configuration, on which the vector fields (or rather,
the $\mathbf{C}^*$-actions which they generate) are
holomorphic, so we can compute the inner product algebro-geometrically
and thereby give a purely algebro-geometric definition of the stability
condition.

\section{Example} \label{sec:ex}

The aim of this section is to work out the stability
criterion in a special case and show how it relates to the existence of
extremal metrics. Let $\Sigma$ be a genus two curve, and $\mathcal{M}$ a line
bundle on it with degree one. The same computation can be
carried out when the genus is greater than two and the line bundle has
degree greater than one. Define $X$ to be the ruled surface
$\mathbf{P}(\mathcal{O}\oplus\mathcal{M})$ over $\Sigma$.
T\o{}nnesen-Friedman~\cite{TF97} constructed a family of extremal metrics on
$X$, which does not exhaust the entire K\"ahler cone. We will show that
$X$ is K-unstable (relative to a maximal torus of automorphisms) for the remaining
polarisations.

Since there are no non-zero holomorphic vector fields on $\Sigma$, a
holomorphic vector field on $X$ must preserve the fibres. Thus, the
holomorphic vector fields on $X$ are given by sections of
$\mathrm{End}_0(\mathcal{O}\oplus\mathcal{M})$. Here, $\mathrm{End}_0$
means endomorphisms with trace zero. The vector field given by the
matrix
\[ \left( \begin{array}{cc}
              -1 & 0 \\
	      0 & 1 
	  \end{array} \right)
\]

\noindent generates a $\mathbf{C}^*$-action $\beta$, and it is up to
scalar multiple the only one that does (see Maruyama~\cite{Mar71} for
proofs). Therefore this must be a multiple of the extremal vector field,
which is then given by $\chi =
\frac{F(\beta)}{\langle\beta,\beta\rangle}\beta$. 

The destabilising test-configuration is an example of deformation to
the normal cone of a subvariety, studied by Ross and Thomas~\cite{RT04},
except we need to take into account the extremal $\mathbf{C}^*$-action
as well. We consider the polarisation $L=C+mS_0$ where $C$ is the
divisor given by a fibre, $S_0$ is the zero section (ie. the image of
$\mathcal{O}\oplus\{0\}$ in $X$) and $m$ is a positive constant. We
denote by $S_\infty$ the infinity section, which as a divisor is just
$S_0-C$. Note that $\beta$ fixes $S_\infty$ and acts on the normal
bundle of $S_\infty$ with 
weight 1. We make no distinction between divisors and their associated
line bundles, and use the multiplicative and additive notations
interchangeably, so for example $L^k=kC+mkS_0$ for an integer $k$.

The deformation to the normal cone of $S_\infty$ is given by the blowup
\[\mathcal{X}:=\widetilde{X\times\mathbf{C}}\xrightarrow{\pi}X\times\mathbf{C}
\]
in the subvariety $S_\infty\times\{0\}$. Denoting the
exceptional divisor by $E$, the line bundle $\mathcal{L}_c=\pi^*L-cE$ is
ample for $c\in(0,\epsilon)$, where $\epsilon$ is the Seshadri constant
of $(S_{\infty},L)$ as in~\cite{RT04}. 
In our case, $\epsilon=m$. Thus we obtain a test-configuration
$(\mathcal{X},\mathcal{L})$ with the $\mathbf{C}^*$ action induced by
$\pi$ from the product of the trivial action on $X$ and the usual
multiplication on $\mathbf{C}$. Denote the restriction of this
$\mathbf{C}^*$-action to the central fibre $(X_0,L_0)$ by $\alpha$.

Since the extremal $\mathbf{C}^*$-action fixes $S_\infty$ we obtain
another action on the test-configuration, induced by $\pi$ from the
product of the extremal $\mathbf{C}^*$-action on $X$ and the trivial
action on $\mathbf{C}$. Let us call the induced action on the central
fibre $\beta$. We wish to calculate $F_\chi(\alpha)$ as defined in
(\ref{eq:modfutaki}). For this we use the following decomposition, with
$t$ being the standard coordinate on $\mathbf{C}$:
\begin{eqnarray*}
H^0(X_0,L_0^k)&=&H^0_X(kL-mkS_\infty)\oplus\bigoplus_{i=1}^{mk-ck}
\frac{H^0_X(kL-(mk-i)S_\infty)}{H^0_X(kL-(mk-i+1)S_\infty)}\oplus\\
&&\bigoplus_{j=1}^{ck}t^j
\frac{H^0_X(kL-(ck-j)S_\infty)}{H^0_X(kL-(ck-j+1)S_\infty)}.
\end{eqnarray*}
According to~\cite{RT04} $\alpha$ acts with weight $-1$ on $t$
that is, it acts with weight $-j$ on the summand of index $j$ above.
Also, $\beta$ acts on
\[\frac{H^0_X(kL-lS_\infty)}{H^0_X(kL-(l+1)S_\infty)}\]
with weight $l$, plus perhaps a constant independent of $l$
which we can neglect, since the matrices are normalized to have trace zero
in the formula for the modified Futaki invariant.
The dimension of this space is $k+l-1$ by the Riemann-Roch theorem.
Writing $A_k, B_k$ for the infinitesimal generators of the actions
$\alpha$ and $\beta$ on $H^0(X_0,L_0^k)$ and $d_k$ for the dimension of
this space, we can now compute
\begin{eqnarray*}
d_k &=& \frac{m^2+2m}{2}k^2 + \frac{2-m}{2}k + O(1),\\
\mathrm{Tr}(A_k) &=& -\frac{c^3+3c^2}{6}k^3 + \frac{c^2-c}{2}k^2 + O(k),\\
\mathrm{Tr}(B_k) &=& \frac{2m^3 + 3m^2}{6}k^3 + \frac{m}{2}k^2 + O(k),\\
\mathrm{Tr}(A_kB_k) &=& -\frac{c^4+2c^3}{12}k^4 + O(k^3),\\
\mathrm{Tr}(B_kB_k) &=& \frac{3m^4+4m^3}{12}k^4 + O(k^3).
\end{eqnarray*}

\noindent Using these, we can compute 
\[ F_\chi(\alpha) = F(\alpha)-\langle\alpha,\chi\rangle = 
F(\alpha)-\frac{\langle\alpha,\beta\rangle}{\langle
\beta,\beta\rangle}F(\beta).\]
We obtain
\[
F_\chi(\alpha)=\frac{c(m-c)(m+2)}{4(m^2+6m+6)}\Big[(2m+2)c^2-(m^2-4m-6)c
+m^2+6m+6\Big]. 
\]

\noindent If $F_\chi(\alpha)\leqslant0$
for a rational 
$c$ between 0 and $m$, then the variety is K-unstable (relative to
a maximal torus of automorphisms).
In~\cite{TF97} (page 23) the condition given for the existence of an
extremal metric of 
the type studied, is that a certain polynomial 
\[ P(\gamma) = \frac{(ka-\gamma)(\gamma-a)}{24}\left[\gamma^2 +
a\left(\frac{-k^2 + 2k + 1}{2k}\right)\gamma + ka^2\right]
\]

\noindent is positive
for $a<\gamma<ka$ in the notation of~\cite{TF97}, where $k$ is the
parameter of the polarization ($k=m+1$ in our notation) and $a$ is a
constant defined in terms of $k$. After a change of variables $k=m+1$
and $\gamma=a(1+c)$ we find that the two conditions are in fact the
same, since a quadratic polynomial with rational coefficients cannot
have an irrational double root.

\bibliographystyle{amsplain} \bibliography{../mybib}

\providecommand{\bysame}{\leavevmode\hbox to3em{\hrulefill}\thinspace}
\providecommand{\MR}{\relax\ifhmode\unskip\space\fi MR }
\providecommand{\MRhref}[2]{%
  \href{http://www.ams.org/mathscinet-getitem?mr=#1}{#2}
}
\providecommand{\href}[2]{#2}
\begin{thebibliography}{10}

\bibitem{Cal82}
E.~Calabi, \emph{Extremal {K}\"ahler metrics}, Seminar on Differential Geometry
  (S.~T. Yau, ed.), Princeton, 1982.

\bibitem{Don97}
S.~K. Donaldson, \emph{Remarks on gauge theory, complex geometry and
  four-manifold topology}, Fields Medallists' Lectures (Atiyah and Iagolnitzer,
  eds.), World Scientific, 1997, pp.~384--403.

\bibitem{Don02}
\bysame, \emph{Scalar curvature and stability of toric varieties}, J.
  Differential Geom. \textbf{62} (2002), 289--349.

\bibitem{DK90}
S.~K. Donaldson and P.~B. Kronheimer, \emph{The geometry of four-manifolds},
  OUP, 1990.

\bibitem{Fut83}
A.~Futaki, \emph{An obstruction to the existence of {E}instein-{K}\"ahler
  metrics}, Invent. Math. \textbf{73} (1983), 437--443.

\bibitem{FM95}
A.~Futaki and T.~Mabuchi, \emph{Bilinear forms and extremal {K}\"ahler vector
  fields associated with {K}\"ahler classes}, Math. Ann. \textbf{301} (1995),
  199--210.

\bibitem{Mab04_1}
T.~Mabuchi, \emph{Stability of extremal {K}\"ahler manifolds}, Osaka J. Math.
  \textbf{41} (2004).

\bibitem{Mar71}
M.~Maruyama, \emph{On automorphism groups of ruled surfaces}, J. Math. Kyoto
  University \textbf{11-1} (1971), 89--112.

\bibitem{MFK94}
D.~Mumford, J.~Fogarty, and F.~Kirwan, \emph{Geometric invariant theory},
  Springer-Verlag, 1994.

\bibitem{RT04}
J.~Ross and R.~P. Thomas, \emph{A study of the {H}ilbert-{M}umford criterion
  for the stability of projective varieties}, preprint (2004).

\bibitem{Tian97}
G.~Tian, \emph{K\"ahler-{E}instein metrics with positive scalar curvature},
  Invent. math \textbf{137} (1997), 1--37.

\bibitem{TF97}
C.~T\o{}nnesen-Friedman, \emph{Extremal {K}\"ahler metrics on ruled surfaces},
  Ph.D. thesis, Odense University, 1997.

\end{thebibliography}

\end{document}